\documentclass[final,a4paper]{amsart}

\usepackage{amsmath,amssymb,amsthm,amsfonts}
\usepackage{epsfig}
\usepackage{enumitem}

\newtheorem{prop}{Proposition}[section]
\newtheorem{theor}{Theorem}[section]

\newtheorem{remark}{Remark}[section]

\newcommand{\n}{\hbox{$\scriptstyle\hbox{\rm
      l\kern-0.22em N}$}}
\newcommand{\N}{\hbox{$ I\kern -0.23em N$}}

\newcommand{\R}{\mathbb R}

\def\Z{\mathbb Z}

\newcommand{\e}{\varepsilon}
\newcommand{\tr}{\hbox{\rm Tr}}

\def\PS #1 #2{\langle #1, #2 \rangle}

\def\YtoXandYinOMEGA{{\displaystyle{\mathop{\scriptstyle{y\to
x}}_{y\in\Omega}}}}

\def\1{1\!\!1}
\def\hyp#1{{\rm\bf (H#1)}}
\def\omegb{{\overline \Omega}}
\def\domeg{{\partial \Omega}}
\def\ou{\overline u}
\def\uu{\underline u}

\def\xb{\bar x}
\def\xh{\hat x}
\def\Fb{\bar F}

\def\SS{{\mathcal S}^N}

\def\pt{\tilde p}
\def\qt{\tilde q}
\def\Xt{\tilde X}
\def\Yt{\tilde Y}

\def\ue{u_{\e}}
\def\xe{x_{\e}}
\def\ze{z_{\e}}

\def\xoe{\frac{x}{\e}}
\def\hyp#1{{\rm\bf (H#1)}}

 1

\begin{document}
\title[{$C^{0,\alpha}$--regularity for superquadratic equations and applications}]{A short proof of the $C^{0,\alpha}$--regularity of viscosity subsolutions for superquadratic viscous Hamilton-Jacobi equations and applications}
\author{Guy Barles}
\address{Laboratoire de Math\'ematiques et Physique Th\'eorique
CNRS UMR 6083, F\'ed\'eration Denis Poisson,
Universit\'e Fran\c{c}ois Rabelais, Parc de Grandmont,
37200 Tours, France} 
\email{barles@lmpt.univ-tours.fr}
\thanks{This work was partially supported by the ANR ``Hamilton-Jacobi et th\'eorie KAM faible'' (ANR-07-BLAN-3-187245)}

\date{}

\maketitle

\begin{abstract}
Recently I.~Capuzzo Dolcetta, F.~Leoni and A.~Porretta obtain a very surprising regularity result for fully nonlinear, superquadratic, elliptic equations by showing that viscosity subsolutions of such equations are locally H\"older continuous, and even globally if the boundary of the domain is regular enough. The aim of this paper is to provide a simplified proof of their results, together with an interpretation of the regularity phenomena, some extensions and various applications.
\end{abstract}

\medskip


\section{Introduction}

In \cite{CDLP},  Capuzzo Dolcetta, Leoni and Porretta obtain a very surprising regularity result for fully nonlinear,  {\em superquadratic}, elliptic equations 
which can be described very easily in their main example which is the one of viscous Hamilton-Jacobi Equations like
\begin{equation}
- \tr (A(x)D^2 u) + \vert Du\vert^p+ \lambda u = f(x) \quad\hbox{in  }\Omega \; ,
    \label{vhje1}
\end{equation}
where $\Omega$ is an open subset of $\R^N$, $\lambda \geq 0$, $p>2$ and $A, f$ are continuous functions taking values respectively in the set of non-negative, $N\times N$ symmetric matrices and $\R$. They show that if $u$ is a locally bounded, upper semicontinuous viscosity subsolution of (\ref{vhje1}) then $u$ is locally H\"older continuous with exponent $\displaystyle \alpha = \frac{p-2}{p-1}$. Furthermore they prove that the local $C^{0,\alpha}$-bound depends only on local $L^\infty$-bounds on $A$, $f$ and $\lambda u^-$. They also provide global regularity results in the case when the boundary $\domeg$ has a sufficient regularity.\\

These results are very unusual and surprising since they provide the regularity of {\em subsolutions} of {\em degenerate}  equations with a {\em superquadratic} growth in $Du$, whereas most of the regularity results for elliptic equations concern {\em solutions} of {\em uniformly elliptic} equations with suitable (subquadratic) {\em growth conditions}. At this point, it is worth mentioning the famous work of Lasry \& Lions \cite{LaLi} where Equation~(\ref{vhje1}) is studied in full details, both in the sub and superquadratic cases, when the second-order term is the laplacian ($A\equiv Id$) : several local gradient bounds are provided by using the Bernstein's method (see also Lions \cite{Llivre} for results in this directions) together with various estimates on the solutions, and all these properties are used to prove existence, uniqueness results in different contexts (infinite boundary conditions, data which are blowing up at the boundary, ergodic problem,...). Most probably, some of their results are still true for (\ref{vhje1}) even if we allow $A$ to degenerate but, at least, their regularity results are valid only for solutions.\\

Coming back to \cite{CDLP}, the way the authors explain it is through the case $A\equiv 0$ for which one has obviously a Lipschitz bound for subsolutions (since $Du$ is clearly bounded) and for general $A$, the viscous Hamilton-Jacobi Equation can be seen as a perturbation of the first-order equation and keeps, at least partially, a similar property when the power of $Du$ is large enough, namely larger than $2$.\\

The aim of the present paper is threefold:\\
(i) to give a slightly simpler proof of this result in a more general setting,\\
(ii) to provide an interpretation of such property in terms of ``state-constraint problems'',\\
(iii) to use this result to obtain, for superquadratic equations, new results for the generalized Dirichlet problem (in the sense of viscosity solutions), for ergodic problems and homogenization problems.\\

In order to be more specific, we come back to Equation (\ref{vhje1}) and we examine again the case $A\equiv 0$ : if $u$ is a subsolution of this equation, then
$$
\vert Du\vert^p \leq  f(x) - \lambda u \quad\hbox{in  }\Omega\; ,$$
and if we assume also that $\lambda =0$, we have a gradient bound which is independent of the $L^\infty$-norm of $u$. And the same property is true for $\lambda \neq 0$ if $u$ is bounded from below.\\

One does not expect such property to be true for elliptic equations and, in general, all the $C^{0,\alpha}$ or Lipschitz bounds depend on (local) $L^\infty$-bounds on $u$. But, as we already mention it above, the authors prove in \cite{CDLP} that the $C^{0,\alpha}$-bound is still true for general $A$ under the same conditions as for the first-order equation.\\

Our approach, whose general framework is described in Section~\ref{GF}, shows why both situations are very similar : roughly speaking, if $u$ is a subsolution of a general equation, we are not going to argue directly on this equation but on a simpler equation for which $u$ is still a subsolution; for the above first-order equation, clearly the only important information is that
$$ \vert Du\vert^p \leq  \vert \vert  f \vert \vert_\infty + \vert \vert  \lambda u^-\vert \vert_\infty\; ,$$
where the $L^\infty$-norm is either a local or a global norm, and this step can be seen as a replacement of a complicated equation by a simpler one. As this (very simple) example shows it, this replacement may depend on (local) $L^\infty$-bounds of $u$ (in the case $\lambda \neq 0$) but once this step is done then the (local) $L^\infty$-bounds will play not role anymore.\\

In order to obtain the $C^{0,\alpha}$-bounds, the key argument consists in building, for the new equation, a family of supersolutions $(w_r)_r$ in balls of radius $r \ll 1$ : these functions are used to control from above the local variations of the subsolution and, of course, this control gives the H\"older regularity. Two points have to be emphasized : first, the $w_r$ are constructed in such a way that they are supersolutions up to the boundary of the balls (this is called ``state-constraints boundary conditions'') and this point is crucial to have a control of the subsolution which is independent of its $L^\infty$-bounds when the new equation does not depend on such $L^\infty$-bounds (the case when $\lambda = 0$  in (\ref{vhje1})). Next the construction of such a family of $w_r$ is possible only in the superquadratic case : we address this question, with several variants, at the end of Section~\ref{GF}.\\

Therefore, in ``good'' cases (typically when $\lambda = 0$ in (\ref{vhje1})), one can obtain $C^{0,\alpha}$-bounds which are independent of any $L^\infty$-bounds on the subsolution and if $\partial \Omega$ is regular enough, these bounds hold up to the boundary. Section~\ref{ME} is devoted to provide various examples of equations to which the framework of Section~\ref{GF} applies and we formulate a general result in Section~\ref{GEN} in which we obtain local modulus of continuity which are not necessarely of H\"older type.\\

Concerning the applications, we are not going to describe them in this introduction; we refer the reader to the corresponding sections. Section~\ref{DP} is devoted to the study of the generalized Dirichlet problem for general superquadratic elliptic equations : assuming or not that the equation is uniformly elliptic, one cannot solve in general the classical Dirichlet problem : we refer for example to Da Lio and the author \cite{BDL2} where the evolution problem is studied and where it is shown that loss of boundary data can occur. For (\ref{vhje1}), it is even obvious that the Dirichlet problem cannot be solved in a classical way since, for smooth enough boundary, the solution is expected to be $C^{0,\alpha}$ up to the boundary and therefore a solution of the classical Dirichlet problem can exist only in cases when the boundary data satisfies rather restrictive conditions. We refer to \cite{CDLP} where this question in studied in full details. On the contrary, we concentrate on solving the generalized Dirichlet problem in the sense of viscosity solutions. The role of $C^{0,\alpha}$-property in this setting is to provide the continuity up to the boundary of the subsolutions which is a key property to obtain comparison results. We refer to \cite{BP1, BP2, BDL1, BRS} for more details.\\

For ergodic and homogenization problems, the role of $C^{0,\alpha}$-bounds is well-known : it is a key argument to solve ergodic problem/cell problem and we show how this can be done for superquadratic equations in sections~\ref{ergo} and ~\ref{ho}. We refer to the bibliography for various
references on ergodic and homogenization problems.

\section{The Key Idea and Main Examples}

\subsection{The General Framework}\label{GF}

The aim of this section is to present a general framework to prove local estimates for viscosity subsolutions of general fully nonlinear elliptic equations with super-quadratic growth. Such equations are written in the form
\begin{equation}\label{GFNLE}
F(D^2 u, Du, u, x) = 0  \quad\hbox{in  }\Omega\; ,
\end{equation}
where $F : \SS \times \R^N \times \R \times \Omega \to \R $ is a continuous function, $\SS$ denoting the space of $N \times N$ symmetric matrices. We assume that $F$ satisfies the (degenerate) ellipticity condition : for any $(p , r , x)\in\R^N \times\R\times \Omega$ and for any $X,Y\in\SS$, 
$$
F(X,p,r,x) \leq F(Y,p,r,x)\quad\hbox{if }X\geq Y.
$$

The function $u$ being a given locally bounded, upper semi-continuous viscosity subsolution  of (\ref{GFNLE}), we make the following assumptions in which $B_r(x)$ denotes the open ball of center $x$ and radius $r$, $d$ is the distance function to $\partial\Omega$ and $a$ is a parameter which can take the value $1$ or $2$.\\

\noindent \hyp{1-$a$} There exists $r_0>0$ and, for any $0 < r \leq r_0$, there exists a continuous function {$G_{r} : \SS \times \R^N \to \R$} satisfying the ellipticity condition such that, for any $x\in \Omega$ with $d (x) \geq ar $,  $G_{r} (D^2 u, Du)\leq 0 $ in $B_r (x)$.\\

In general, the functions $G_r$ depend on local or global $L^\infty$-bounds of $u$ and this is the role of the parameter ``$a$'' to express in which way : if $a=2$, then, for any $x\in \Omega$ such that $d (x) \geq 2r $,  $B_r (x) \subset \Omega_r : = \{y:\ d(y)>r\}$ and a priori $G_r$ depends on the $L^\infty$-bounds of $u$ on $\Omega_r$. On the contrary, if $a=1$, then the balls $B_r (x)$ cover the whole domain $\Omega$ and either $G_r$ does not depend on any $L^\infty$-bounds of $u$ (typically if $F$ does not depend on $u$) or on a global $L^\infty$-bounds of $u$.\\

The next assumption is\\

\noindent \hyp{2}  For any $0<r\leq r_0$, there exists $w_r \in C(\overline{B_r(0)})$ such that $w_r(0) = 0$,  $w_r \geq 0$ in $B_r(0)$ and
\begin{equation}\label{SCwr}
 G_r (D^2 w_r, D w_r) \geq \eta_r > 0 \quad\hbox{on  }\overline{B_r(0)}\setminus \{0\}\; ,
\end{equation}
for some $\eta_r >0$.\\

We notice that $w_r$ can be assumed to be radially symmetric since one can replaced it by $\tilde  w_r$ defined for $0\leq s \leq r$ by
$$\tilde  w_r (s) := \inf_{ \vert y \vert =s}w_r(y)\; .$$
By standard arguments, it is easy to show that $\tilde  w_r$ still satisfies (\ref{SCwr}) because it is essentially an infimum of supersolutions.
In \noindent \hyp{2}, it is also important to remark that $w_r$ (or $\tilde  w_r$) is a supersolution up to the boundary (state-constraint boundary condition) and this point plays a key role to obtain bounds which are independent of $L^\infty$-bounds of $u$ when $G_r$ has this property.\\

The last assumption is \\

\noindent \hyp{3} If $v$ is a bounded upper semicontinuous viscosity subsolution of $G_r (D^2 v, D v) \leq 0$ in  $B_r(0)\setminus \{0\}$, then $v(y)\leq v(0) + w_r (y)$ in $B_r(0)$.\\

In other words, Assumption \hyp{3} means that, for any $r$, one has a comparison result for the state-constraint problem in $B_r(0)\setminus \{0\}$ and that we can compare any upper semicontinuous viscosity subsolution with the strict viscosity supersolution $w_r$.\\

We give more comments on \hyp{2}-\hyp{3} at the end of this section.\\

We have the
\begin{prop}\label{basic} If \hyp{1-2}-\hyp{2}-\hyp{3} hold then the viscosity subsolution $u$ of (\ref{GFNLE}) satisfies the following property : for any $x \in \Omega$ and any $y \in B_r (x)$ where $r \leq d(x)/2$, we have
$$ u(y) \leq u(x) + w_r (y-x)\; .$$
In particular, for any $\delta >0$, $u$ is uniformly continuous on $\overline{\Omega_{\delta}}$ and $w_{\delta/2}$ is a modulus of continuity of $u$ on $\overline{\Omega_\delta}$ (recall that $w_{\delta/2}(s)$ depends only on $ \vert s \vert $).\\
Finally, if $\Omega$ is a $C^{1,1}$-domain and if $w_r(s)\leq \tilde K s^\alpha$ for some $\tilde K>0$, $0<\alpha \leq 1$ independent of $r$, then $u$ can be extended as a $C^{0,\alpha}$-function on $\overline\Omega$.
\end{prop}

By $C^{1,1}$-domain, we mean that the distance function $d$ is $C^{1,1}$ in a neighborhood of $\partial\Omega$, say in $\{x:\ d(x) < \delta_0\}$, and therefore $n(x)=-Dd(x)$ is Lipschitz continuous in this neighborhood.\\

\proof \ The beginning of the proof is obvious since $y \mapsto u(x +y)$ is a viscosity subsolution of the $G_r$-equation in $B_r(0)\setminus \{0\}$ and, by \hyp{3}, we have the first part of the result.

Next, if $x,y \in \overline{\Omega_\delta}$ satisfy $ \vert x-y \vert <\delta/2$, then we have at the same time $y \in B_{\delta/2} (x)$, $d(x)\geq 2(\delta/2)$ and $x \in B_{\delta/2} (y)$, $d(y)\geq 2(\delta/2)$ ; hence in the above result we can exchange the roles of $x$ and $y$ and get
$$ \vert u(y) - u(x) \vert \leq w_{\delta/2} (y-x)\; ,$$
which proves the uniform continuity statement.

Finally, if we assume that $\Omega$ is a $C^{1,1}$-domain and $w_r(s)=\tilde Ks^\alpha$ for some $\tilde K>0$, $0<\alpha \leq 1$ independent of $r$, we first 
estimate $\vert u(x) - u(x-\delta n(x)) \vert$ for $x \in \Omega$ such that $d(x) \leq \delta_0/2$ and $\delta \leq \delta_0/2$, $\delta_0$ being defined as above as a constant such that $d$ is $C^{1,1}$ on $\{z:\ d(z) < \delta_0\}$.

To do so, we introduce the points defined, for $k\in \N$, by
$$x_k = x-\frac{\delta}{2^k}n(x)\; ;$$
hence $x_0=x-\delta n(x)$ and $\lim_{k\to +\infty}\,x_k = x$. We estimate
$$ u(x_K)-u(x-\delta n(x)) =u(x_K) - u(x_0)= \sum_{k=1}^{K} [u(x_k)- u(x_{k-1})]\; .$$
Since $\displaystyle \vert x_k - x_{k-1} \vert =  \frac{\delta}{2^k}$, $\displaystyle d(x_{k-1}) = d(x) + \frac{\delta}{2^{k-1}} > \frac{ \delta}{2^k}$ and $\displaystyle d(x_k) = d(x) + \frac{\delta}{2^k} > \frac{ \delta}{2^k}$, we have
$$
 \vert u(x_k)- u(x_{k-1})  \vert \leq \tilde K \vert x_k - x_{k-1} \vert^\alpha =  \tilde K\frac{\delta^\alpha}{(2^\alpha)^k}\; .
$$
Therefore
$$ \vert  u(x_K)- u(x)\vert  \leq \tilde K  \delta^\alpha\sum_{k=1}^{K}\frac{1}{(2^\alpha)^k} \leq \bar K \delta^\alpha\; ,$$
and letting $K$ tends to $+\infty$ shows that $\vert u(x)- u(x-\delta n(x)) \vert \leq \bar K \delta^\alpha$ since the sum is converging.

Next we consider $x,y \in \Omega$ such that $ \vert  x-y \vert < \delta_0/4$ and we want to estimate  $\vert u(y) - u(x) \vert$. If either $d(x) \geq \delta_0/2$ or $d(y) \geq \delta_0/2$, this can be done by the second step of the proof above since $x,y \in \overline{\Omega_{\delta_0/4}}$, hence by $\tilde K \vert  x-y \vert^\alpha$. Therefore the interesting case is when
$d(x) < \delta_0/4$ and $d(y)< \delta_0/4$. In this case, we argue in the following way~: we introduce $0 < \delta \leq \delta_0/2$ and write
$$ u(y) - u(x) = [u(y)-u(y-\delta n(y))] + [u(y-\delta n(y))-u(x-\delta n(x))] + [u(x-\delta n(x))-u(x)] \; .$$
Denoting by $m_\delta$ the modulus of continuity of $u$ in $\overline{\Omega_\delta}$ and using the above result to estimate the first and third term, we obtain
$$ \vert u(y) - u(x) \vert \leq  2\bar K \delta^\alpha +m_\delta( \vert (y-\delta n(y))-(x-\delta n(x))\vert) \; .$$
And by the regularity of $\Omega$, $\vert (y-\delta n(y))-(x-\delta n(x))\vert \leq (1+L\delta) \vert  x-y \vert \leq (1+L\delta_0) \vert  x-y \vert$, $L$ being the Lipschitz constant of $n$ in $\Omega_{\delta_0}$. Recalling that $m_\delta(t) = \tilde K t^\alpha$ for $t\leq \delta$, we see that the choice $ \delta = (1+L\delta_0) \vert  x-y \vert$ provides the answer (notice that changing $\delta_0$ into a smaller constant, we can assume, without loss of generality, that $1+L\delta_0 \leq 2$). And the proof is complete.$\Box$\\

Before checking \hyp{2}-\hyp{3} in various cases, we provide further comments on them. If $F$ is uniformly elliptic then, in general, the $G_r$ are also uniformly elliptic and, in most cases, it is impossible to build such functions which are (viscosity) supersolution up to the boundary : more precisely, if the $G_r$ are subquadratic, the results of Da Lio\cite{DL} imply that the equation cannot hold up to the boundary and therefore we have no hope to construct the $w_r$ in this case. On the contrary, we can indeed have a state-constraint boundary condition in the superquadratic case but in a strange way : in fact, as it is noticed in \cite{LaLi}, this state-constraint boundary condition is equivalent to $\displaystyle \frac{\partial {w_r}}{\partial n} = + \infty$ on $\partial B_r(0)$, which implies that we do not have any supersolution requirement on the boundary since $w_r - \phi$ cannot achieve a minimum on $\partial B_r(0)$ if $\phi$ is a smooth function.\\

As a consequence, the reader who is not very familiar with viscosity solutions theory but wants to be convinced by the results (at least for (\ref{vhje1}) or for uniformly elliptic equations and for smooth solutions), can check the following assumption instead of \hyp{2}

\noindent \hyp{2'}  For any $0<r\leq r_0$, there exists $w_r \in C(\overline{B_r(0)})$ such that $w_r(0) = 0$,  $w_r \geq 0$ in $B_r(0)$ and
\begin{equation}\label{SCwrS}
 G_r (D^2 w_r, D w_r) \geq \eta_r > 0 \quad\hbox{in  }B_r(0)\setminus \{0\}\; ,
\end{equation}
\begin{equation}\label{SCwrB}
\displaystyle \frac{\partial {w_r}}{\partial n} = + \infty \quad\hbox{on  }\partial B_r(0)\; ,
\end{equation}
for some $\eta_r >0$.\\

Once you have checked such property by building smooth functions $w_r$ inside $B_r(0)\setminus \{0\}$ and assuming that we only consider smooth subsolutions (to prove, for example, that one has uniform $C^{0,\alpha}$-bounds), the checking of \hyp{3} is immediate since, if $v$ is (at least) Lipschitz continuous on $\overline{B_r(0)}$, $v-w_r$ cannot achieved its maximum on the boundary because of (\ref{SCwrB}) and standard Maximum Principle type arguments allows to conclude that the maximum can be achieved only at $0$, which provides the desired result.\\

Finally, one may wonder if the same approach could work in the subquadratic case replacing, in \hyp{2}, the boundary condition on $w_r$ by the natural condition
$$ w_r(y) \to + \infty  \quad\hbox{when }|y| \to r\; .$$
In fact, it is hopeless to obtain similar results as shown by the following example : we consider, for $N\geq 3$ and $1\leq p \leq 2$, the equation
$$
-\Delta u + |Du|^p = 0 \quad\hbox{in  }B_1(0)\; .
$$
It is easy to check that, if $C\alpha \leq 1$ and $\alpha$ is small enough, the functions $C|x|^\alpha$ are subsolutions of this equation. In particular, if the functions $w_r$ could be built, we would have a uniform control on the local modulus of continuity of these subsolutions : this is clearly not the case when considering the particular sequence $\alpha^{-1}|x|^\alpha$ as $\alpha \to 0$.

\subsection{The Main Example}\label{ME}

In this section, we are going to prove that the following type of viscous Hamilton-Jacobi Equation
\begin{equation}
- \tr (a(x)D^2 u) + H(x,Du)+ c(x) u = f(x) \quad\hbox{in  }\Omega\; ,
    \label{vhje}
\end{equation}
enters into the framework we described in the previous section provided that\\

{\bf (i)} The function $x \mapsto a(x)$ is a continuous function defined on $\overline \Omega$, with values in the space of $N \times N$ matrices, such that 
$$ a(x)p \cdot p \geq 0 \quad\hbox{ for any   }p\in \R^N\; .$$

{\bf (ii)} The function $(x,p) \mapsto H(x,p)$ is a continuous function defined on $\overline \Omega \times \R^N$ and there exists constants $K_1, K_2$ and $m>2$ such that 
$$H(x,p) \geq K_1 \vert p \vert ^m -K_2\; ,$$
for any $x \in \overline \Omega$ and $ p \in \R^N$.

{\bf (iii)} The functions $x \mapsto c(x), f(x)$ are real-valued, continuous functions defined on $\overline \Omega$.\\

Our result is the
\begin{theor}\label{regu} If {\bf (i)-(iii)} holds then Assumptions \hyp{1-2}-\hyp{2}-\hyp{3} are satisfied and $w_r(s) \leq \tilde K_r s^\alpha$ for for $\displaystyle \alpha = \frac{m-2}{m-1}$ and some constant $\tilde  K_r >0$. As a consequence, any locally bounded subsolution $u$ of (\ref{vhje}) is in $C^{0,\alpha}_{loc}(\Omega)$. Moreover, if $\Omega$ is a $C^{1,1}$-domain and if $u$ is bounded on $\overline \Omega$, then $u \in C^{0,\alpha}(\overline \Omega)$.\\
Finally, if the function $x \mapsto c(x)u(x)$ is bounded from below on $\overline \Omega$, then any locally bounded subsolution of (\ref{vhje}) is globally bounded and therefore in $C^{0,\alpha}(\overline \Omega)$ if $\Omega$ is a $C^{1,1}$-domain.
\end{theor}

The assumption ``the function $x \mapsto c(x)u(x)$ is bounded from below'' may seem strange; in fact, it is satisfied in two interesting particular cases : the first one is when $c \equiv 0$ and the second one is when $c(x)\geq 0$ and $u$ is bounded from below (in particular $u\geq 0$ on $\overline \Omega$).\\

\proof \  If $u$ is a subsolution of (\ref{vhje}), it is a subsolution of
$$ - \vert\vert a \vert\vert_\infty \sum_{\lambda_i(D^2u) > 0}\lambda_i(D^2u)  + K_1 \vert Du \vert ^m \leq R\; ,$$
for some $R$ large enough. Therefore the functions $G_r(p,M)$ are of the form
$$ G_r(M,p):= - \vert\vert a \vert\vert_\infty \sum_{\lambda_i(M) > 0}\lambda_i(M)  + K_1 \vert p \vert ^m - R
\; ,$$
for some $R$ large enough, to be specified later and since they are all of the same form (just $R$ may change), we are going to use the simplified notation $G$ only.

In order to build the functions $w_r$, we first build $w_1$. To do so, for $C_1, C_2>0$ to be chosen later on and for $\alpha = \frac{m-2}{m-1}$, we consider the function
$$ w_1 (x):= C_1  \vert x \vert ^\alpha + C_2 (d^\alpha(0) - d^\alpha(x))\; ,$$
where $d$ is equal to the distance to $\partial B_1 (0)$ (i.e. $d(x)=1- \vert x \vert $) if, say, $ \vert x \vert \geq 1/2$ and we regularize it in $B_{1/2}(0)$ by changing it into $\varphi(1- \vert x \vert )$ where $\varphi$ is a smooth, non-decreasing and convex function such that $\varphi (s) $ is constant for $s \leq 1/4$ and $\varphi(s)=s$ for $s\geq 1/2$. With this change, the new function $d$ is smooth in $\overline{B_1(0)}$.

Obviously we have $w_1 (0) = 0$, $w_1 \geq 0$ in $\overline{B_1(0)}$ and $w_1$ is smooth in $\overline{B_1(0)}\setminus \{0\}$, which will allow us to prove that $w_1$ is a supersolution of $G \geq 0$ in $B_1(0)\setminus \{0\}$ by just computing derivatives. For the boundary of the ball, we have $\displaystyle \frac{\partial w_1}{\partial n} = +\infty$ and it is immediate that, if $\phi$ is a smooth function, $w_1-\phi$ cannot achieve a local minimum on $\partial B_1(0)$ and the viscosity supersolution property holds since there is no constraint.

In $B_1(0)\setminus \{0\}$, we compute the derivatives of $w_1$
$$ Dw_1 (x) = \alpha C_1  \vert x \vert ^{\alpha -2} x - \alpha C_2 d^{\alpha-1}(x) Dd(x)\; ,$$
\begin{eqnarray*}
D^2w_1 (x) & = & \alpha C_1  \vert x \vert ^{\alpha -2}Id +  \alpha (\alpha-2) C_1  \vert x \vert ^{\alpha -4} x\otimes x \\
 & & - \alpha C_2 d^{\alpha-1}(x) D^2 d(x)- \alpha(\alpha -1) C_2 d^{\alpha-2}(x) Dd(x)\otimes Dd(x)\; .
\end{eqnarray*}
To simplify the computation, we make several remarks : on one hand, $-Dd(x) = \mu(x) x$ for some $\mu(x)\geq 0$; this is a consequence of the way we built the function $d$. As a consequence, we have
\begin{eqnarray*}
\vert Dw_1 (x)\vert ^m & = & (\vert \alpha C_1  \vert x \vert ^{\alpha -2} x  \vert +\vert \alpha C_2  d^{\alpha-1}(x) Dd(x) \vert    )^m \\
& \geq  & \vert \alpha C_1  \vert x \vert ^{\alpha -2} x  \vert ^m
+ \vert \alpha C_2  d^{\alpha-1}(x) Dd(x) \vert ^m\; .
\end{eqnarray*}
On the other hand, using that $d(x)= \varphi (1- \vert x \vert )$, with $\varphi$ convex, we obtain
\begin{eqnarray*}
D^2w_1 (x) &\leq &\alpha C_1  \vert x \vert ^{\alpha -2}Id +\alpha C_2 d^{\alpha-1}(x) \left( \varphi '' + 
\frac{\varphi'}{ \vert x \vert }\right) \frac{x}{ \vert x \vert } \otimes \frac{x}{ \vert x \vert }  \\
& & + \alpha(1-\alpha) C_2 d^{\alpha-2}(x) Dd(x)\otimes Dd(x)\; ,
\end{eqnarray*}
and 
$$\lambda_i (D^2w_1 (x)) \leq \alpha C_1  \vert x \vert ^{\alpha -2} +\alpha C_2 d^{\alpha-1}(x) \left( \varphi '' + 
\frac{\varphi'}{ \vert x \vert }\right)+ \alpha(1-\alpha) C_2 d^{\alpha-2}(x) \vert Dd(x)\vert^2\; .$$

These properties imply that, in order to prove the expected inequality for $w_1$, we can (almost) consider the two terms separately. More precisely, taking into account the value of $\alpha$, the $C_1  \vert x \vert ^\alpha$-term yields
$$ -  \alpha  \vert\vert a \vert\vert_\infty C_1  \vert x \vert ^{\alpha -2} + K_1  \vert \alpha C_1  \vert x \vert ^{\alpha -2} x  \vert ^m =  \vert x \vert ^{\alpha -2} \left(-  \alpha  \vert\vert a \vert\vert_\infty C_1 + K_1   \alpha^m C_1 ^m\right)\; .$$
By choosing $C_1$ large enough, this quantity can be as large as we wish on $B_1(0)\setminus \{0\}$. On the other hand, the $C_2 (d^\alpha(0) - d^\alpha(x))$-term yields
$$ -  C_2 \vert\vert a \vert\vert_\infty d^{\alpha-2}(x) \left(
\alpha d(x) \left( \varphi '' + 
\frac{\varphi'}{ \vert x \vert }\right)+ \alpha(1-\alpha)  \vert Dd(x)\vert^2 \right)+ K_1  \vert \alpha C_2  d^{\alpha-1}(x) Dd(x) \vert ^m \; .$$
Here we have to consider two cases : either $ \vert x \vert Ê\geq 1/2$ and then $\varphi' =1$, $\varphi'' =0$ and $Dd(x) = - x/ \vert x \vert $; therefore the above quantity is nothing but
$$ -  C_2 \vert\vert a \vert\vert_\infty d^{\alpha-2}(x) \left(
\alpha\frac{d(x)}{ \vert x \vert }+ \alpha(1-\alpha) \right)+ K_1 \alpha^m C_2^m  d^{m(\alpha-1)}(x)\; .$$
Since $m(\alpha-1)=\alpha-2$, by choosing $C_2$ large enough, this quantity can be positive (and even greater than $\tilde k d^{\alpha-2}(x)$ for any $\tilde k >0$). Finally, for $ \vert x \vert \leq 1/2$, the above quantity is bounded.

In order to conclude for $w_1$, the above computations shows that, by taking first $C_2$ large enough and then $C_1$ large enough, then $G (D^2w_1,Dw_1) \geq 1$ on $\overline{B_1(0)}\setminus \{0\}$, where $1$ can be replaced by any positive constant.

Next step consists in building $w_r$ using $w_1$. To do so, we set 
$$ w_r(x):= r^\alpha w_1(\frac{x}{r})\; .$$
It is easy to check that for $0<r \leq 1$, $G ( D^2  w_r , D  w_r) \geq 1$ on $\overline{B_r(0)}\setminus \{0\}$. In fact, the ``$1$'' as well as the ``$R$'' in $G$ can be replaced by $r^{\alpha -2} \geq 1$ and $r^{\alpha -2}R \geq R$ respectively.

To conclude the checking of \hyp{1-2}-\hyp{2}, we remark that a subsolution of (\ref{vhje}) is a subsolution of (say) $G=0$ for $R= \vert\vert (f-cu)^+ \vert\vert_\infty$ where, for the local estimates, the $L^\infty$-norm is taken on balls of the form $B_r(x)$ for $x \in \Omega$ and $r\leq d(x)$ where $d$ denotes here the distance of $x$ to $\partial \Omega$.

It remains to check \hyp{3}. There are several way to do it. Taking the above construction of $w_r$ into account, the simplest one consists probably in using the arguments of Da Lio and the author\cite{BDL1} :  instead to comparing directly the upper semicontinuous viscosity subsolution $v$ of $G = 0$ in  $B_r(0)\setminus \{0\}$ with $v(0) + w_r (y)$, one compares it with $v(0) + w_{r'} (y)$ for some $r' < r$; this can be done by considering 
$$ \max_{(x,y) \in \overline{B_{r}(0)} \times \overline{B_{r'}(0)}}\,\left(v(x)- ( v(0) + w_{r'} (y)) - \frac{ \vert x-y \vert ^2}{\e^2}\right)\; .$$
Choosing $\e \ll r-r'$ and, using the fact that $v$ and $v(0) + w_{r'}$ are bounded, it is easy to show that the maximum is achieved for $x \in B_{r}(0)$ (and not on the boundary of the ball) and for $y\in \overline{B_{r'}(0)}$. Using in an essential way the fact that $v(0) + w_{r'} (y)$ is a (strict) supersolution up to the boundary allows to perform all the usual viscosity solutions arguments in this context. We refer to \cite{BDL1} for details.

Once  \hyp{1-2}-\hyp{2}-\hyp{3} hold, we can apply Proposition~\ref{basic} to obtain all the $C^{0,\alpha}$-properties on $u$ : indeed, it suffices to remark that $w_1$ is $C^{0,\alpha}$ and therefore one has $w_1(x) \leq \tilde K  \vert x \vert ^\alpha$ and, by construction, $w_r (x) \leq \tilde K  \vert x \vert ^\alpha$. A priori $\tilde K$ depends on $R$, hence on the ball we consider inside $\Omega$.

The last point comes from the fact that the $C^{0,\alpha}$-bounds on $u$ depend on $R=\vert\vert (f-cu)^+ \vert\vert_\infty$ : if $cu$ is bounded from below, so is $R$ and the H\"older bound is independent of the ball. It is then easy to show that $u$ is globally bounded as soon as it is locally bounded.$\Box$\\

We conclude this section by describing further examples of applications of the above result. We first consider the pde
\begin{equation}\label{qle}
 - \frac{\Delta u}{1+ \vert Du \vert ^k} +  \vert Du \vert ^m + c(x) u = f(x) \quad\hbox{in  }\Omega\; ,
\end{equation}
where $c$ and $f$ satisfy the same assumptions as above and $k,m>0$.

Multiplying the equation by $1+ \vert Du \vert ^p$, we (almost) recover the same framework as above with $a = Id$ and
$$ H(x,u,p)= (1+ \vert p \vert ^k)  \vert p \vert ^m +  (1+ \vert p \vert ^k) c(x) u - f(x)(1+ \vert p \vert ^k)\; .$$
Clearly, by using some local $L^\infty$-bounds on $u$, such $H$ satisfies an assumption like {\bf (ii)} if $k+m >2$. And Theorem~\ref{regu} applies almost readily. This example shows that one may have to rewrite the equation before applying the above result.

\subsection{A general result}\label{GEN}

We have presented only semilinear or quasilinear examples but it is clear that Theorem~\ref{regu} extends to the case of fully nonlinear equations (\ref{GFNLE}) since only the functions $G_r$ are playing a role and not the function $F$ itself. On the other hand, we can treat the case of more general nonlinearities than just superquadratic power of $Du$.

In order to formulate our general result, we introduce the class $\mathcal{P}$ of functions $ h : [0,+\infty) \to [0,\infty)$ satisfying : \\
(i) $h$ is a $C^1$, convex function,\\
(ii) $\displaystyle t \mapsto \frac{h(t)}{t^2}$ is non-decreasing for $t\geq 1$,\\
(iii) $\displaystyle \int_1^{+\infty} \frac{t}{h(t)}dt < +\infty$.\\

We first notice that all the functions $t \mapsto t^m$ for $m>2$ are in the class $\mathcal{P}$. We also use below the fact that any function $h \in \mathcal{P}$ satisfies\\
(iv)   $\displaystyle \int_1^{+\infty} \frac{t^2h'(t)}{[h(t)]^2}dt < +\infty$\\
a result which is obtained by a simple integration by part, remarking that (ii) implies that the function $\displaystyle t \mapsto \frac{h(t)}{t^2}$ has a limit when $t\to +\infty$ and that, by (iii), this limit is necessarely $+\infty$.
 
 Our result is the

\begin{theor}\label{mainregu} Assume that the following assumption holds :\\

\noindent \hyp{4} for any $R>0$,
$$ F(M,p,u,x)\leq 0\quad  \Rightarrow \quad \tilde G_R (M) + K^R_1 h(\vert p \vert) -K^R_2 \leq 0 \; ,$$
for any $M \in \SS$, $ p \in \R^N$, $ \vert u \vert Ê\leq R$, $x \in \overline \Omega$, where $K^R_1,K^R_2$ are positive constant, $h$ is a function in the class $\mathcal{P}$ and $\tilde G_R$ is a Lipschitz continuous function satisfying the ellipticity condition, which is homogeneous of degree $1$.\\
Then, any locally bounded subsolution $u$ of (\ref{GFNLE}) is uniformly continuous in $\overline{\Omega_\delta}$ with a modulus of continuity depending only on $\delta$, $R_\delta:=\vert \vert u \vert\vert_{L^\infty(\overline{\Omega_{\delta/2}})}$ and the different constants and functions appearing in \hyp{4} with $R=R_\delta$. Moreover, if $h(t)=t^m$, $\Omega$ is a $C^{1,1}$-domain and either $u$ is bounded on $\overline \Omega$ or $K^R_1,K^R_2$ can be choosen independent of $R$, then $u \in C^{0,\alpha}(\overline \Omega)$.\\
Finally, if $h(t)=t^m$ and either $u$ is bounded on $\overline \Omega$ or if $K^R_1,K^R_2$ can be chosen independent of $R$, then any locally bounded subsolution of (\ref{GFNLE}) is globally bounded and therefore in $C^{0,\alpha}(\overline \Omega)$ if $\Omega$ is a $C^{1,1}$-domain.
\end{theor}

Theorem~\ref{mainregu} has the most general formulation we can provide but, of course, on each particular case, it may be interesting to look more precisely at the right inequality to be used on $F$ and which may depend on the subsolution.\\

\proof \  It follows very closely the proof of Theorem~\ref{regu} : if $u$ is a subsolution of (\ref{GFNLE}), it is also (locally) a subsolution of
$$  \tilde G_R ( D^2u)  + K^R_1 h(\vert Du \vert) \leq K^R_2 \; ,$$
for some $R$ large enough. This inequality defines our functions $G_r$ that we just denotes by $G$ as above.

In order to build the functions $w_r$, we only build $w_1$ and we will indicate later on how to argue for $w_r$; it is worth pointing out that, for general functions $h$, the scaling argument of the proof of Theorem~\ref{regu} cannot be used.

We introduce the $C^2$, increasing functions $\chi_1, \chi_2 : [0,1] \to \R$ defined by $\chi_1 (0)=0, \chi_2(0)=0$, $\chi'_1 (0)=+\infty, \chi'_2(0)=+\infty$ and
$$ h(\chi'_1 (t)) = \frac{\chi'_1 (t)}{t}\quad , \quad \chi''_2 (t) = - h(\chi'_2 (t))\; .$$
We refer to the appendix where the existence of $\chi_1, \chi_2$ is studied and various properties are obtained that we are going to use in the computations below.

For $C_1, C_2>0$ to be chosen later on, we consider the function
$$ w_1 (x):= C_1 \chi_1 ( \vert x \vert ) + C_2 (\chi_2(0) - \chi_2(d(x))\; ,$$
where $d$ is as in the proof of Theorem~\ref{regu}.

Obviously we have $w_1 (0) = 0$, $w_1 \geq 0$ in $\overline{B_1(0)}$ and $w_1$ is smooth in $\overline{B_1(0)}\setminus \{0\}$, which will allow us to prove that $w_1$ is a supersolution of $G\geq 0$ in $B_1(0)\setminus \{0\}$ by just computing derivatives. For the boundary of the ball, we have $\displaystyle \frac{\partial w_1}{\partial n} = +\infty$ because $\chi'_2(0)=+\infty$.

We next remark that by the convexity of $h$, $h(t_1 + t_2) \geq h(t_1) + h(t_2)$ for any $t_1,t_2 \geq 0$ because the function $s \to h(s + t_2) - h(s)$ is an non-decreasing function on $[0,+\infty)$ and therefore it achieves its minimum at $s=0$. This allows us to treat separately the $\chi_1$ and $\chi_2$ terms.

For the $\chi_1$-term, denoting by $\xh$ the quantity $ \vert x \vert ^{-1}x$, the second derivative is 
$$ C_1  \chi'_1 ( \vert x \vert ) D^2( \vert x \vert ) + C_1 \chi''_1 ( \vert x \vert ) \xh \otimes \xh\; .$$
But, $\chi_1$ is concave (see the Appendix) and taking into account, the Lipschitz continuity of $\tilde G_R$, the contribution of this $\chi_1$-term in $G$ is estimated by
$$ - C_1 L_R \frac{\chi'_1 ( \vert x \vert )}{ \vert x \vert } +  K^R_1 h(C_1\chi'_1 ( \vert x \vert ))\; ,$$
where $L_R$ is the Lipschitz constant of $ \tilde G_R $.
But, by the second property of the class $\mathcal{P}$, $h(C_1\chi'_1 ( \vert x \vert )) \geq C_1^2 h(\chi'_1 ( \vert x \vert ))$, and we are left with:
$$ (K^R_1 C_1^2 - C_1 L_R) h(\chi'_1 ( \vert x \vert )) \; $$
and since $h(\chi'_1 ( \vert x \vert ))$ is bounded away from $0$, this term can be as large as we wish by choosing $C_1$ large.

Now we turn to the $\chi_2$-term and we first examine the case $ \vert x \vert  \geq 1/2$ where $ \vert Dd(x) \vert =1$. The $G$ quantity is estimated by
$$ - C_2 L_R \chi'_2 (d(x)) D^2d(x) + C_2 L_R \chi''_2 (d(x)) + K^R_1 h(C_2 \chi'_2 (d(x)))\; ,$$
and by similar arguments as for the $\chi_1$-term, we can transform it into
$$ - C_2 L_R M_2 \chi'_2 (d(x)) + C_2 L_R \chi''_2 (d(x)) + C_2^2 K^R_1 h(\chi'_2 (d(x)))\; ,$$
where $M_2$ stands for the $L^\infty$-norm of $D^2 d$ for $1/2 \leq  \vert x \vert  \leq 1$. Using the superquadratic behavior of $h$, it is clear that this quantity is positive by choosing $C_2$ large enough. And as in the proof of Theorem~\ref{regu}, $d$ is smooth in $B_{1/2}(0)$ and therefore this quantity is bounded from below.

The conclusion follows for $w_1$, as in the proof of Theorem~\ref{regu}, by first taking $C_2$ large enough and then $C_1$ large enough.

Next, in order to build $w_r$, we set 
$$ w_r(x):= r w_1(\frac{x}{r})\; .$$
By using the homogeneity of $\tilde G_R$, it is clear that we have just to repeat the above construction of $w_1$ with $K^R_1, K^R_2$ being replaced by $ rK^R_1,r K^R_2$.

\begin{remark}
It is worth pointing out that, if $h$ satisfies $h(ct)\geq \rho(c)h(t)$ for any $c,t \geq 1$ for some function $\rho$ such that $c^{-2}\rho(c) \to +\infty$ as $c \to +\infty$, then the constants $C^r_1,C^r_2$ associated with $rK^R_1,r K^R_2$ satisfy $r C^r_1,rC^r_2 \to 0$ as $r \to 0$ and $w_r(x) \to 0$ uniformly on $\overline{B_r(0)}$ as $r \to 0$. We obtain in this case a very good control of the subsolution, even when we are close to the boundary of $\domeg$, as in the $h(t)=t^m$-case (which is a particular case in which this condition holds). Unfortunately it is not very difficult to check that this condition implies that, actually, $h(t) \geq t^m$ for some $m>2$, at least for large $t$, and therefore we are in the $h(t)=t^m$-case.
\end{remark}

To conclude the checking of \hyp{1-2}-\hyp{2}, we remark that a subsolution of (\ref{GFNLE}) is a subsolution of the above equation in a ball $B_r(x)\subset\subset \Omega$ for $R= \vert\vert u \vert\vert_{L^\infty(B_r(x))}$ and, for \hyp{3}, the arguments of Da Lio and the author\cite{BDL1} still apply.

Finally the case $h(t) = t^m$ is treated exactly with the argument of Theorem~\ref{regu}, which apply readily.$\Box$\\

We conclude this section technical remarks which appear already in \cite{CDLP}. In Section~\ref{GF}, we present a general framework in which balls play a central role; but, if $\Omega$ is a convex, $C^3$  domain, the proof of Theorem~\ref{regu} (and, in the same way, of Theorem~\ref{mainregu}) can be done directly in $\Omega$ by showing that, for any $x,y \in \omegb$ 
$$ u(y) \leq u(x) +  C_1  \vert y-x \vert ^\alpha + C_2 (d^\alpha(x) - d^\alpha(y))\; .$$
To do so, the key point is that, for fixed $x$, the right-hand side is a supersolution of the $G$-equation with a state-constraint boundary condition of the boundary of $\Omega$ and the proof of this fact just use the inequality $(y-x) \cdot Dd(y) \leq 0$ on $\omegb$ which is true in convex domains (and that remains true after a suitable regularization of $d$).

An other variant which is more local and can be useful close to the boundary of $\Omega$ consists in using rectangles instead of balls, in particular when the subsolution $u$ is known to be bounded : we write the unit rectangle as $\{x=(x',x_N)\, ; \vert x' \vert < 1 \; , \; -1 < x_N < 1\}$ where $x'$ denotes $(x_1, \cdots,x_{N-1})$. In the above construction, one may replace the distance to the boundary by $\min(1-x_N, x_N+1)$ which is just the distance to the parts of the boundary $\{x_N = 1\}$ and $\{x_N = - 1\}$ (of course, one has to regularize it). Built in that way (with exactly the same arguments), $w_1$ (and the $w_r$) are supersolutions but up to the boundary only for the parts $\{x_N = 1\}$ and $\{x_N = - 1\}$; for the other part of the boundary, namely $\{\vert x' \vert = 1\}$, one has to manage in order to have $u(x) + w_r(y-x) \geq u(y)$ and this is where we use the $L^\infty$-bound on $u$; this inequality can be obtained without any difficulty by taking $C_1$ large enough. Close to the boundary --which can be flatten if it is smooth enough-- such argument allows to obtain a more precise behavior of the subsolution.

\section{Application 1 : The Dirichlet Problem for Superquadratic Elliptic Equations}\label{DP}

We consider in this section the Dirichlet problem consisting in solving (\ref{GFNLE}) together with the boundary condition
\begin{equation}
    u(x) = g(x) \quad \hbox{on  }\domeg \; ,
    \label{GDBC}
\end{equation}
where $g$ is a continuous function.\\

To simplify matter, we are going to assume that $F$ satisfies \hyp{4} with $h(t) = t^m$ for some $m>2$ and we call \hyp{4'} this new assumption.\\

For different reasons, and here the regularity of $g$ will be an unusual addition reason, it is well-known that, in general, (\ref{GFNLE})-(\ref{GDBC}) has no solution assuming the boundary data in a classical sense and one has to use the formulation of the generalized Dirichlet boundary condition in the viscosity solutions sense which reads
\begin{equation}\label{minbc}
    \min(F(D^2 u, Du, u, x) , u-g)\leq 0~~\mbox{on 
    $\partial\Omega,$}
\end{equation}
and
\begin{equation}\label{maxbc}
    \max(F(D^2 u, Du, u, x) , u-g)\geq 0~~\mbox{on
    $\partial\Omega.$}
\end{equation}
Roughly speaking, these relaxed conditions mean that the equations has
to hold up to the boundary, when the boundary condition is not assumed
in the classical sense.  In general, the key argument to justifies them is that
they appear naturally when one passes to the limit in the vanishing
viscosity method using typically the well-known ``half-relaxed limits method''~; 
here one may think more on truncation arguments on the superquadratic dependence in $Du$. 

For superquadratic, {\em uniformly elliptic} equations, these boundary conditions reduces to $u-g \leq 0$ on $\partial\Omega$ as we will see it below and the fact that the viscosity solution inequality (\ref{maxbc}) is ``inactive'' in the sense that $u-\phi$ cannot have a minimum point on $\partial\Omega$ for any smooth function $\phi$; formally this typically means that $\displaystyle \frac{\partial u}{\partial n} = + \infty$. We refer to Da Lio and the author\cite{BDL2} for more comments in this direction.

We first briefly analyze the loss of boundary conditions for (\ref{GFNLE})-(\ref{GDBC}).

From now on, we assume that $\Omega$ is a smooth, bounded domain with a $C^{3}$-boundary. We denote by $d$
a $C^{3}$-function agreeing in a neighborhood  $\mathcal{W}$ of $\partial\Omega$ with
the signed distance function to $\partial\Omega$ which is positive in $\Omega$ and negative in
$\R^N\setminus\overline\Omega$ and we denote by $n(x):=-Dd(x)$ for
all $x\in\mathcal{W}.$ If $x\in\partial\Omega,$ $n(x)$ is just the
unit outward normal to $\partial\Omega$ at $x.$

\begin{prop}\label{sub}
Assume that $F$ is a continuous function satisfying the ellipticity condition and \hyp{4'}, and that $g\in C(\partial \Omega)$. If $u$ is a bounded, usc subsolution of (\ref{GFNLE})-(\ref{GDBC}), then
$$ u \leq \varphi \quad \hbox{on  }\domeg \; .$$
\end{prop}

Under the assumption of Proposition~\ref{sub}, it is clear that (\ref{GFNLE})-(\ref{GDBC}) have no solution (even viscosity solution) assuming the boundary data continuously : indeed, because of \hyp{4'}, the subsolutions (and therefore the solutions) are expected to be H\"older continuous up to the boundary and if $g$ is not in the right $C^{0,\alpha}$-space, then $u$ cannot be equal to $g$ on the boundary.

Since we assume $u$ to be only usc, artificial discontinuities may appear on the boundary : indeed if $u(x) < g(x)$ at some point $x\in\partial\Omega$, then one may replace the value of $u(x)$ by any value between $u(x)$ and $g(x)$ : in that way, the function remains usc and still satisfies the boundary condition. To avoid this difficulty, we always assume that, for any $x\in\partial\Omega$
\begin{equation}\label{sub-bound}
u(x) =
\limsup_{\YtoXandYinOMEGA}\,u(y)\: .
\end{equation}

\smallskip

\proof\  We use a result of Da Lio\cite{DL} : if $u(x_0)>g(x_0)$ at some point $x_0\in\partial\Omega$, then one has 
\begin{equation}\label{subp1}
    \begin{array}{c}
    \displaystyle{ \liminf_{\displaystyle{\mathop{\scriptstyle{y\to
    x_{0}}}_{\alpha\downarrow 0}}}}
\left\{          \displaystyle
\left[F\left(
-{1\over{\alpha^2}} ( Dd(y)\otimes Dd(y) + o(1)),
{{Dd (y) + o(1)}\over{\alpha}},
u(y),y
\right) \right]\right\}
            \leq 0
                \end{array}
\end{equation}

But, because of \hyp{4'} which we use with $R= \vert\vert  u \vert\vert_\infty$, the $F$-quantity is larger than
$$\tilde G_R (-{1\over{\alpha^2}} ( Dd(y)\otimes Dd(y) + o(1)) ) + K^R_1 \vert {{Dd (y) + o(1)}\over{\alpha}} \vert ^m -K^R_2\; ,$$
which clearly tends to $+\infty$ since $\tilde G_R $ is homogeneous of degree $1$, $K^R_1>0$ and $m>2$.$\Box$\\

To state a comparison result for the generalized Dirichlet problem, we introduce the following structure assumptions on $F$ which are inspired from Rouy, Souganidis and the author \cite{BRS}~: this is natural since we need here a comparison result between {\em continuous}, and even $C^{0,\alpha}$, subsolutions and (a priori) discontinuous supersolutions. Of course, the framework of \cite{BRS} has to be adapted to take into account, at the same time, the regularity of subsolutions and the superquadratic growth of the equation in $Du$. In these assumptions, $\alpha = \frac{m-2}{m-1}$.\\

The assumptions are the following

\medskip \noindent \hyp{5} \hspace{0.3cm} For all $R$ there exists
$\gamma_R>0$ such that, for all $x \in \omegb $, $-R \leq v\leq u \leq
R$, $p\in\R^N$ and $M\in \SS$

\smallskip
$$ F(M,p,u,x)-F(M,p,v,x) \geq \gamma_R
(u-v)\;  .$$

\bigskip
\noindent \hyp{6}\hspace{0.3cm} There exists $0< \beta <\alpha $ and for all $0< \mu < 1$, $R, K >0$, there exists a
function $m_{R,K}:\R^+ \to \R$ such that $m _{R,K} (t)\to 0$ when $t \to 0$, and a constant $C(\mu)$ such that, for all $\varepsilon >0$

\smallskip
$$ F(Y,q,u,y)- \mu F(\mu^{-1} X,\mu^{-1} p,\mu^{-1}u,x)\leq
m_{R,K}\left(
(1-\mu) +C(\mu)\varepsilon
\right)$$

\vspace{0.1cm}
\noindent
for all $x,y \in\omegb$, $\vert  u \vert  \leq R$,
$p,q \in \R^N$ and for all matrices $X, Y \in \SS$
satisfying the following properties
\begin{equation}\label{inegmatthbis}
-K\varepsilon^{-2 +\alpha - 2\beta } Id
\leq \left(\begin{array}{cc} X & 0 \\ 0 & -Y
\end{array}\right) \leq K \varepsilon^{-2 +\alpha - 2\beta}
\left(\begin{array}{cc} Id
& -Id \\ -Id & Id
\end{array}\right) +  K\varepsilon^{\alpha - 2\beta} Id ,
\end{equation}
\begin{equation}
	\vert p-q \vert \leq K \varepsilon (1+ \vert p \vert +\vert q
\vert)\; , \; \vert p \vert +\vert q
\vert \leq K\varepsilon^{\alpha-1-\beta}\; ,
	\label{propgrad}
\end{equation}
\begin{equation}
		\vert x-y \vert \leq K \varepsilon .
	\label{propxy}
\end{equation}

 \bigskip 
 
Our result is the following
\begin{theor}\label{compres}{\bf : }Assume that \hyp{4'}-\hyp{6} hold and let $u$ and $v$ be respectively a bounded usc
subsolution and a bounded lsc super-solution of (\ref{GFNLE})-(\ref{GDBC}), $u$ satisfying condition (\ref{sub-bound}). Then
$$ u \leq v \quad\hbox{ on } \omegb\; .$$
\end{theor}

We refer to \cite{BRS} for comments on the different assumptions. We recall here that we adapt them to take into account the facts that $u$ is H\"older continuous and $F$ is superquadratic in $p$. We illustrate these assumptions on the equation
$$
- \tr (a(x)D^2 u) + b(x) \vert Du \vert ^m + c(x) u = f(x) \quad\hbox{in  }\Omega\; ,
$$
where $a,b,c$ and $f$ are continuous functions with $a = \sigma \sigma^T$ for some continuous function $\sigma$, where $\sigma^T$ denotes the transpose matrix of $\sigma$. Of course, \hyp{5} means that $c(x) >0$ on $\overline{\Omega}$. Next, for \hyp{6}, we have to estimate
$$ Q:=F(Y,q,u,y)- \mu F(\mu^{-1} X,\mu^{-1} p,\mu^{-1}u,x) $$
for $F(X, p,u,x) =- \tr (a(x)X) + b(x) \vert p \vert ^m + c(x) u - f(x)$. We have
$$ Q = - \tr (a(y)Y-a(x)X)) + (b(y) \vert q \vert ^m - \mu b(x) \vert  \mu^{-1} p \vert ^m)+ (c(x)-c(y)) u - (f(y)-\mu f(x))\; .$$
To simplify the checking, we remark that we can estimate each term separately : indeed, if we have a $m_{R,K}^1$ modulus for one term with a $C^1(\mu)$ constant and a $m_{R,K}^2$ modulus for an other term with a $C^2(\mu)$ constant, then, $m_{R,K}^1,m_{R,K}^2$ can be choosen as being increasing functions, the sum of these two terms satisfies the assumption with $m_{R,K}^1+m_{R,K}^2$ and $C^1(\mu)+C^2(\mu)$. Finally, for the dependence in $\beta$, we point out that, if a term satisfies \hyp{6} for some $\beta$, it satisfies it for any $\beta' < \beta$; hence we can take the smallest $\beta$ which appears for the different term.

For the first, third and fourth, \hyp{6} follows from standard arguments, namely
$$ \vert \tr (a(y)Y-a(x)X))\vert \leq K \varepsilon^{-2 +\alpha - 2\beta}\vert \sigma (x) - \sigma(y)\vert^2 + 4K\vert\vert \sigma\vert\vert_\infty \varepsilon^{\alpha - 2\beta}\; ,$$
and the assumption holds if $\sigma \in C^{0,\gamma}(\overline{\Omega})$ with $\gamma > 1 - \alpha/2$ by taking $\beta$ small enough,
$$\vert c(x)-c(y)) u \vert + \vert f(y)-\mu f(x)\vert \leq m_c( \vert x-y \vert)R + m_f( \vert x-y \vert) + \vert 1-\mu  \vert. \vert \vert f \vert \vert_\infty\; ,$$
where $m_c,m_f$ are some modulus of continuity for $c$ and $f$ respectively. Here we have no constraint for \hyp{6} to be satisfied.

We end up with the unusual term $ b(y) \vert q \vert ^m-\mu b(x) \vert  \mu^{-1} p \vert ^m$. We set $\tilde \mu = \mu^{1-m}>1$. We have
$$
b(y) \vert q \vert ^m - {\tilde \mu} b(x) \vert p \vert ^m  =  (b(y)-b(x)) \vert q \vert ^m + b(x)(  \vert q \vert ^m - {\tilde \mu}  \vert p \vert ^m)
$$
but, by (\ref{propgrad}), 
$$\vert q \vert \leq (1+K \varepsilon) (1-K \varepsilon)^{-1}\vert p \vert + K\varepsilon (1-K \varepsilon)^{-1}$$
and therefore, if $\e$ is small enough compared to $\tilde  \mu - 1$
$$
b(y) \vert q \vert ^m - {\tilde \mu} b(x) \vert p \vert ^m  \leq  (b(y)-b(x)) \vert q \vert ^m + \vert\vert b \vert\vert_\infty \frac{( 1- {\tilde \mu})}{2}  \vert p \vert ^m + O(\e^m)\; ,
$$
and if $b\in C(\omegb)$, the right-hand side of this inequality is estimated by $O(\e^m)$, for $\e$ small enough.\\

\begin{remark} In order to check \hyp{6}, it may be more convenient to change $F$ by multiplying it by a positive quantity. For example, in the case of (\ref{qle}), multiplying the equation by $1+ \vert Du \vert ^k$ leads to an equation which still satisfies \hyp{5} (and even a stronger property) and for which the checking of \hyp{6} is easier. Two remarks on that example : (i) in order to formulate an optimal result (\ref{qle}), it would be necessary to take into account the stronger version of \hyp{5}; we are not going to do it in order to avoid long and technical details. And (ii) to treat (\ref{qle}) we need an additional argument to check \hyp{6}. The above argument on the $b(x) \vert p \vert ^m$ leading term, shows that it is estimated by $\vert\vert b \vert\vert_\infty \frac{( 1- {\tilde \mu})}{2}  \vert p \vert ^m + O(\e^m)$ and therefore we have a ``good'' term which can be used to estimate other terms (recall that $\tilde \mu >1$). For example, and again we choose a simple example for the sake of clarity, a linear term $d(x)\cdot p$ can be estimated in the following way, using (\ref{propgrad}) and Young's inequality
\begin{eqnarray*}
d(x)\cdot p - d(y)\cdot q & \leq &\vert d(x)-d(y) \vert.\vert p \vert + \vert d(y) \vert .\vert p-q \vert \\
& \leq &\vert d(x)-d(y) \vert.\vert p \vert + \vert d(y) \vert .K\e (1 +  \vert p \vert  +  \vert q \vert )\\
& \leq & C(\mu) \vert d(x)-d(y) \vert^{m'} + \vert\vert b \vert\vert_\infty  \frac{ \vert 1- {\tilde \mu} \vert}{4}  \vert p \vert ^m + o_\e (1)\; ,
\end{eqnarray*}
where $(m')^{-1} + m^{-1} = 1$. Therefore the continuity of the coefficients is enough for all types of terms.
\end{remark}

\medskip

\proof \ We follow readily the proof of \cite{BRS} taking into account the results of Section~\ref{ME} which implies that $u$ is in $C^{0,\alpha}$ (at least with the correct redefinition on the boundary) and Proposition~\ref{sub}. For $0 < \mu < 1$, close to $1$, we consider $M_\mu = \max_{\ \omegb }\,(\mu u-v)$ and  argue by contradiction assuming that $\liminf_{\mu \to 1}\, M_\mu > 0$. 

By the regularity of the boundary, there exists a $C^2$-function $\chi : \R^N \to \R^N$ which is equal to $n$ is a neighborhood of $\partial \Omega$ and
 we introduce the test-function $\Phi_\varepsilon : \omegb \times
\omegb \to \R$ by
$$
\Phi_\varepsilon (z,w) = \mu u(z)-v(w) - \vert \frac{z-w}{\varepsilon}
+\chi(\frac{z+w}{2} )
\vert^k\; , $$
where $k$ is a large power, namely $k \geq \alpha/\beta$.

Let $(x,y)$ be a global maximum point of $\Phi_\varepsilon$ on $\omegb \times
\omegb$.  For notational simplicity here we drop the dependence of $x$
and $y$ on $\varepsilon$.

By standard arguments, since $u$ is $C^{0,\alpha}$, we have
$$ \vert \frac{x-y}{\varepsilon}
+\chi(\frac{x+y}{2} )
\vert^k \leq C \vert x-y\vert ^\alpha\; ,$$
and since $\displaystyle \frac{x-y}{\varepsilon}$ is bounded, this means that
$$ \vert \frac{x-y}{\varepsilon}
+\chi(\frac{x+y}{2} )
\vert^k \leq \tilde C \e^\alpha\; .$$
Using this estimate, tedious but straightforward computations shows that the elements $(p,X) \in \overline J^{2,+} \mu u(x)$ and $(q,Y) \in \overline J^{2,-}
v(y)$ given by the Jensen-Ishii's Lemma, satisfy all the properties listed in \hyp{6} and
$$ F(Y,q,v(y),y)\geq 0\quad \hbox{and}\quad \mu F(\mu^{-1} X,\mu^{-1} p,\mu^{-1}u(x) ,x)\leq 0\; .$$
An immediate application of \hyp{5}-\hyp{6} concludes.$\Box$

Theorem~\ref{compres} leads us to the
\begin{theor}\label{exi-uni}{\bf : }Assume that \hyp{4'}-\hyp{6} hold and that there exists some constant $M>0$ such that
$$F(0,0,-M,x) \leq 0 \leq F(0,0,M,x)\quad\hbox{on  }\omegb\; ,$$
then, for any $g\in C(\domeg)$,  there exists a unique solution $u$ of (\ref{GFNLE})-(\ref{GDBC}); moreover $u \in 
C^{0,\alpha}(\omegb)$. In addition, if $\gamma_R$ can be chosen independent of $R$ in \hyp{5},  there exists a unique solution $u_\infty\in 
C^{0,\alpha}(\omegb)$ of the state-constraints problem and for any viscosity subsolution $w$ of (\ref{GFNLE}) satisfying condition (\ref{sub-bound}), we have
$$ w \leq u_\infty \quad\hbox{ on } \omegb\; .$$
\end{theor}

\proof\ Using \hyp{5}, we can assume that the constant $M$ is larger than $\vert\vert g\vert\vert_\infty$ since we can change it into $\max(M,\vert\vert g\vert\vert_\infty)$. Therefore $-M$ and $M$ are respectively sub and supersolution of (\ref{GFNLE})-(\ref{GDBC}) and one can apply the Perron's method of Ishii \cite{I3} (See also the user's guide for viscosity solutions \cite{cil}) with the version up to the boundary of Da Lio \cite{DL}. The regularity of the solution comes from Theorem~\ref{mainregu}.

The existence of $u_\infty$ is a little bit more complicated : $-M$ (without any change on $M$) is a subsolution. For the supersolution, we need a supersolution  up to the boundary and to do so, we borrow arguments from Tabet Tchamba~\cite{T} : this supersolution, denoted by $\ou$, is built in the following way
$$
\ou(x) = M + K (1+d^\alpha(x))\quad \hbox{on  } \omegb\; ,
$$
for some large constant $K>0$. Using \hyp{5} shows
$$
F(D^2 \ou, D\ou, \ou, x) \geq F(D^2 \ou, D\ou, M, x)
$$
and we can use \hyp{4'} and the ideas of the proof of Theorem~\ref{regu} to show that for $K$ large enough, $\ou$ is a supersolution (up to the boundary) in a neighborhood of $\domeg$. Then we use the fact that $\gamma_R$ can be chosen independent of $R$ in \hyp{5} and the regularity of $d^\alpha(x)$ in $\Omega$ to extend this supersolution property to the whole domain. The rest of the proof is standard.$\Box$\\

\section{Application 2 : The Ergodic Problem}\label{ergo}

The problem is formulated in the following way : find a constant $c$ and a solution $w_\infty \in C(\omegb)$ of the equation
\begin{equation}\label{GFNLE-erg}
F(D^2 w_\infty, Dw_\infty, x) = c  \quad\hbox{in  }\Omega\; ,
\end{equation}
associated to a state-constraint boundary condition on $\domeg$, namely
\begin{equation}\label{GFNLE-erg-sc}
F(D^2 w_\infty, Dw_\infty, x) \geq c  \quad\hbox{on  }\domeg\; .
\end{equation}

We first point out that results in this direction were first obtained for the case of (\ref{vhje1}) in \cite{LaLi}. 
Our result is the
\begin{theor}\label{erg-pb}{\bf : }Assume that $\Omega$ has a $C^3$-boundary and that \hyp{4'} and \hyp{6} hold. There exists a unique constant $c$ such that the state-constraints problem (\ref{GFNLE-erg}) has a continuous viscosity solution $w_\infty$. Moreover, if $(F_k)_k$ is a sequence of continuous functions, satisfying the ellipticity condition and \hyp{4'}-\hyp{6} uniformly with respect to $k$, and if $F_k \to F$ locally uniformly, then the associated ergodic constants $c_k$ converge to $c$.
\end{theor}

\proof\ We just sketch it since it follows from very standard arguments. By Theorem~\ref{exi-uni}, for any $\lambda >0$, there exists a unique solution $w_\lambda \in C^{0,\alpha}(\omegb)$ of the state-constraints problem for the equation
\begin{equation}\label{GFNLE-erg-l}
F(D^2 w_\lambda, Dw_\lambda, x) + \lambda w_\lambda= 0  \quad\hbox{in  }\Omega\; .
\end{equation}
Moreover, since $-\lambda^{-1}\vert\vert F(0,0,x) \vert \vert_\infty$ is a subsolution of this problem, we have $w_\lambda \geq -\lambda^{-1}\vert\vert F(0,0,x) \vert \vert_\infty$ in $\omegb$ and it follows that
$$
F(D^2 w_\lambda, Dw_\lambda, x)  \leq \vert\vert F(0,0,x) \vert \vert_\infty  \quad\hbox{in  }\Omega\; .
$$
By Theorem~\ref{mainregu}, we have uniform $C^{0,\alpha}$-estimates for the functions $w_\lambda$ (uniform in $\lambda$) and if $x_0$ is any point of $\omegb$,
$ \tilde w_\lambda (x) := w_\lambda (x)-w_\lambda(x_0)$ is also uniformly bounded.

Finally, examining carefully the proof of Theorem~\ref{exi-uni} and the construction of $\ou$, it follows also that $\lambda w_\lambda$ is actually bounded and not only bounded from below. It follows that, by Ascoli's Theorem, we can extract a uniformly converging subsequence from $(\tilde w_\lambda)_\lambda$ and, since $(- \lambda w_\lambda (x_0))_\lambda$ is bounded, we can also assume that it converges along the same subsequence. Denoting by $w_\infty$ and $c$, the limits of $(\tilde w_\lambda)_\lambda$ and $(- \lambda w_\lambda (x_0))_\lambda$ respectively and remarking that $ \tilde w_\lambda$ solves
$$
F(D^2  \tilde w_\lambda, D \tilde w_\lambda, x) +\lambda  \tilde w_\lambda = - \lambda w_\lambda (x_0)  \quad\hbox{in  }\Omega\; ,
$$
we can pass to the limit and this shows that $w_\infty$ and $c$ solve the ergodic problem (\ref{GFNLE-erg})-(\ref{GFNLE-erg-sc}).\\

The uniqueness of $c$ comes from standard arguments comparing solutions $(w_\infty,c)$ and $(\tilde w_\infty, \tilde c)$ : the arguments of the proof of Theorem~\ref{compres} shows that, if we use the comparison arguments with $w_\infty$ playing the role of the subsolution and $\tilde w_\infty$ playing the role of the supersolution, then $\tilde c \leq c$; and uniqueness comes from the symmetric roles of $w_\infty$ and $\tilde w_\infty$.

For the last part, if $w^k$ is a solution of the ergodic problem $F_k=c_k$ with $w^k (x_0)=0$, then the $w^k$ satisfy uniform $C^{0,\alpha}$ bounds and the $c_k$ are also bounded. By Ascoli's Theorem, we can extract a uniformly converging subsequence $(w^{k'})_{k'}$ from $(w^k)_k$ and we can also assume that the associated $(c_{k'})_{k'}$ is converging. By the stability result for viscosity solution, the limit is a solution of the $F$-ergodic problem, showing that the limit of $(c_{k'})_{k'}$ is necessarely $c$, by the uniqueness of the ergodic constant. This proves that the limit of any subsequence of 
$(c_{k})_{k}$ is converging tp $c$ and therefore $c$ is the limit of the bounded sequence $(c_{k})_{k}$.$\Box$

\begin{remark} Of course the key argument of the proof of Theorem~\ref{erg-pb} is the $C^{0,\alpha}$ bound on the $w_\lambda$ and it is clear that
the same result holds with different type of boundary conditions, either Neumann boundary conditions or in the periodic setting.
\end{remark}

An immediate consequence of Theorem~\ref{erg-pb} concerns the large time behavior of solutions of evolution equations 

\begin{theor}\label{ltb}{\bf : }Assume that $\Omega$ is a bounded domain with a $C^3$-boundary and that \hyp{4'} and \hyp{6} hold. Then, for any $u_0\in C(\omegb)$, there exists a unique viscosity solution of the Cauchy-Dirichlet problem
\begin{equation}
     u_{t} + F(D^2 u, Du, x)
    = 0 \quad\hbox{in  }\Omega \times (0,+\infty)\; ,
    \label{vhje-evol}
\end{equation}
\begin{equation}
    u(x,0) = u_{0}(x) \quad \hbox{on  }\domeg\; ,
    \label{id-evol}
\end{equation}
\begin{equation}
    u(x,t) = u_{0}(x) \quad \hbox{on  }\domeg \times (0, + \infty)\; ,
    \label{bc-evol}
\end{equation}
and we have
$$ \lim_{t \to \infty}\, \frac{u(x,t)}{t} = - c^+\; ,$$
uniformly on $\omegb$, where $c$ is the constant given by Theorem~\ref{erg-pb}. Moreover, if $c<0$ and $F(M,p,x)$ is convex in $(M,p)$, then $u(x,t)$ converges locally uniformly in $\Omega$ to the unique solution of the stationary Dirichlet problem
\begin{equation}
F(D^2 w, Dw, x)
    = 0 \quad\hbox{in  }\Omega\; ,
    \label{vhje-stat}
\end{equation}
\begin{equation}
 w(x) = u_{0}(x) \quad \hbox{on  }\domeg\; .
 \label{bc-stat}
\end{equation}
\end{theor}

\proof.\ We just sketch the proof since it is tedious but easy adaptation of the argument of \cite{BDL2}, \cite{BRS} and \cite{T}. The existence and uniqueness of $u$ comes directly from the argument of \cite{BDL2}.\\

If $c >0$, one remarks that, if $w_\infty$ is a solution of the ergodic problem (\ref{GFNLE-erg})-(\ref{GFNLE-erg-sc}) and if $K:=  \vert\vert u_0 -w_\infty \vert\vert_{L^\infty (\omegb)}$, then $-ct +w_\infty+ K$ is a supersolution of (\ref{vhje-evol})-(\ref{bc-evol}) which is above $u_0$ at $t=0$ and therefore, by comparison (using, for example, the arguments of \cite{BRS} recalled in Section~\ref{DP}), we have, for all $t>0$
$$ u(x,t) \leq -ct+w_\infty(x) + K\quad\hbox{on  }\omegb \times (0,+\infty)\; .$$
We point out that this argument is based on the fact that (\ref{GFNLE-erg-sc}) holds and, of course, (\ref{bc-evol}) is understood in the viscosity sense.

On the other hand,  $-ct +w_\infty - K$ is a subsolution of (\ref{vhje-evol})-(\ref{bc-evol}) which is below $u_0$ at $t=0$ and on $\domeg$, and applying again a comparison result, we have
$$ -ct+w_\infty(x) - K \leq u(x,t) \quad\hbox{on  }\omegb \times (0,+\infty)\; .$$
 and the result follows from the two above inequalities for $c>0$.\\
 
If $c \leq 0$, then we can use $w_\infty - K $ as subsolution, while if, $\xb$ is a point far enough to $\Omega$, we can take functions like $\vert x- \xb \vert^2$ as supersolutions because of \hyp{4'} (see \cite{T}). It follows
$$ w_\infty(x) - K \leq u(x,t) \leq \vert x- \xb \vert^2 \quad\hbox{on  }\omegb \times (0,+\infty)\; ,$$
and therefore $u$ is uniformly bounded and, obviously
$$ \lim_{t \to \infty}\, \frac{u(x,t)}{t} = 0\; .$$

If $c<0$, applying the half-relaxed limit method, we obtain that the functions 
$$\ou (x) = \limsup_{{\displaystyle{\mathop{\scriptstyle{y\to
x}}_{t \to +\infty}}}} u (y,t) \quad \hbox{  and  } \quad \uu (x) = \liminf_{{\displaystyle{\mathop{\scriptstyle{y\to
x}}_{t \to +\infty}}}} u (y,t)\; ,$$
are respectively viscosity sub and supersolution of the stationary Dirichlet problem (\ref{vhje-stat})-(\ref{bc-stat}). But $F(M,p,x)$ is convex in $(M,p)$ and $w_\infty$ is a strict viscosity subsolution of the equation, therefore we have a Strong Comparison Result for this Dirichlet problem. As a consequence, we have  
$$ \ou \leq \uu \quad\hbox{in  }\Omega\; .$$
Moreover, by standard arguments, the function $w:=\ou=\uu$ is continuous in $\Omega$ and can be extended as a continuous function up to the boundary which is the unique viscosity solution of (\ref{vhje-stat})-(\ref{bc-stat}). Finally the fact that $\ou=\uu$ in $\Omega$ implies that $u(x,t)$ converges locally uniformly to $w$ in $\Omega$ and the proof is complete.
$\Box$\\

Additional arguments are required to study the behavior of the bounded function $ u(x,t) + ct$ for $c\geq 0$:  we refer to \cite{T} where it is proved that this function converges uniformly to $w_\infty + C$ for some constant $C$ in the case of the equation
$$
     u_{t} -\Delta u + \vert Du \vert^m
    = f(x) \quad\hbox{in  }\Omega \times (0,+\infty)\; ,
$$
where $m>2$ and $f \in C(\omegb)$. In fact, while the result of Theorem~\ref{ltb} only requires Theorem~\ref{erg-pb} and (more or less standard)  comparison results for either the Cauchy-Dirichlet problem ot the stationary Dirichlet problem, the $c\geq 0$-case uses the Strong Maximum Principle (cf. Bardi and Da Lio \cite{bdl2,bdl3} and Da Lio \cite{DL-smp}).

\section{Application 3 : Homogenization Problems}\label{ho}

We are interested in homogenization-singular perturbation problems of the form
\begin{equation}\label{GFNLE-hom}
F(\e D^2 \ue, D \ue, \xoe, x) + \ue = 0  \quad\hbox{in  }\Omega\; ,
\end{equation}
with, say, the Dirichlet boundary condition (in the generalized snse)
\begin{equation}\label{DBC-hom}
 \ue = g  \quad\hbox{on  }\domeg\; .
\end{equation}
We assume that $F(M,p,y,x)$ is a continuous function which satisfies the ellipticity property and is $\Z^N$-periodic in $y$, and
$$ F_\e (M,p,x) = F(M,p,\xoe ,x)\; ,$$
satisfies \hyp{4'} uniformly wrt $\e$ and \hyp{6}.

To study (\ref{GFNLE-hom})-(\ref{DBC-hom}), we have first to introduce the {\em cell problem} which consists in finding, for any $(x,p) \in \omegb \times \R^N$, a $\Z^N$-periodic function $y\mapsto u_1 (p,x,y)$ and a constant $\Fb (p,x)$ which solve
\begin{equation}\label{CP}
F(D_{yy}^2 u_1, D_y u_1 + p , y, x) = \Fb (p,x) \quad \hbox{in  }\R^N\; .
\end{equation}
The existence of such $u_1$ and $\Fb (p,x)$, and the uniqueness of $\Fb (p,x)$ is ensured by the argument of Section~\ref{ergo}.

We first give a result on this cell problem.
\begin{prop} Under the above assumptions on $F$, $\Fb$ is a continuous, coercive function. Moreover, if $(F_k)_k$ is a sequence of continuous functions satisfying the ellipticity condition and the same assumptions as $F$, uniformly with respect to $k$, and converging locally uniformly to $F$, then the associated sequence $(\Fb_k)_k$ also converges locally uniformly to $\Fb$.
\end{prop}

\proof\ This result is mainly a consequence of Theorem~\ref{erg-pb}, and in particular of the last part of Theorem~\ref{erg-pb} which say that the ergodic constant depends continuously of the nonlinearity under suitable conditions. Only the coercivity requirement does not come from this result.

If $u_1$ is a continuous, $\Z^N$-periodic solution of (\ref{CP}), we can consider the $\max_{y\in \R^N} ( u_1(p,x,y) )$ which is achived at some point $\bar y$. Using that $u_1$ is a subsolution of (\ref{CP}), we obtain
$$ F(0, p , \bar y, x) \leq \Fb (p,x) \; ,$$
and the coercivity property of $\Fb$ is an immediate consequence of \hyp{4'}.$\Box$\\

Our result for the homogenization problem is the

\begin{theor}\label{hom} 
Under the above assumptions on $F$, the family $(\ue)_{\e>0}$ converges uniformly on $\omegb$ to the unique solution $\ou$ of
\begin{equation}\label{limeqn}
 \left\{
 \begin{array}{rcl}
\Fb (D\ou,x) + \ou = 0 \quad\hbox{in  } \Omega \, ,\\
\noalign{\vskip6pt}
\ou = g \quad \hbox{on  } \domeg\; .
  \end{array}
 \right.
\end{equation}
\end{theor}

\proof\ We first notice notice that the sequence $(\ue)_{\e>0}$ is in a compact subset of $C(\omegb)$ : indeed, by standard arguments, it is easy to see that
$$ \vert\vert \ue \vert\vert_\infty \leq  \vert\vert F(0,0,y,x) \vert\vert_\infty\; ,$$
and because $F_\e$ satisfies \hyp{4'} uniformly wrt $\e$, we have a uniform $C^{0,\alpha}$-bound.

Up a subsequence, we can therefore assume that $(\ue)_{\e>0}$ converges uniformly on $\omegb$ to a function $\ou$ and we have to prove that $\ou$ solves (\ref{limeqn}). If it is the case, we are done since, by results of Perthame and the author\cite{BP1,BP2}, (\ref{limeqn}) has a unique solution and therefore the whole sequence $(\ue)_{\e>0}$ converges to the unique solution of (\ref{limeqn}).

We just prove that $\ou$ is a subsolution of (\ref{limeqn}), the supersolution property being obtained by similar arguments. Let $\phi$ be a smooth function on $\omegb$ and assume that $\xb$ is a strict maximum of $\ou-\phi$.

If $\xb \in \domeg$ and $\ou (\xb) \leq g(\xb)$, we have nothing to prove. Otherwise we have to show that $\Fb(D\phi(\xb), \xb)+\ou(\xb) \leq 0$. To do so, we use the perturbed test-function method of Evans \cite{LCE1,LCE2} with an additional trick which is already used in Da Lio, Lions, Souganidis and the author \cite{BDLLS} and in \cite{GB}.

We argue in the following way : for $k \gg 1$, we introduce
$$ F_k(M,p,y,x): = \min \{F(M',p',y',x')\ ;\ \vert M'-M \vert+\vert p'-p \vert + \vert y'-y\vert+\vert x'-x\vert \leq k^{-1}\}\; ,$$
and we solve (\ref{CP}) for $F_k$, with $p=D\phi(\xb)$ and $x= \xb$; the solution is denoted by $v_k$ and the right-hand side of (\ref{CP}) by $\Fb_k (D\phi(\xb),\xb)$.

Next we consider maximum points $(x_\e, z_\e)$ of
$$ \ue (x) - \phi(x) - \e v_k (\e^{-1}z) - \frac{\vert x-z\vert^2}{\beta^2}\; ,$$
where $0 < \beta \ll \e$. By standard arguments, $\xe \to \xb$ as $\e \to 0$, and even if $\xe \in \domeg$, $\ue (\xe) > g(\xe)$ by the uniform convergence of $\ue$ to $\ou$. Therefore the Jensen-Ishii's Lemma implies the existence of $(\pt ,\Xt ) \in \overline J^{2,+} \ue (\xe)$ and $(\qt,\Yt) \in \overline J^{2,-}
v_k(\frac{\ze}{\e})$ such that
$$ \left(\begin{array}{cc} \Xt & 0 \\ 0 & -\frac{1}{\e}Y
\end{array}\right) \leq \frac{1}{\beta^2}
\left(\begin{array}{cc} Id
& -Id \\ -Id & Id
\end{array}\right) +  \left(\begin{array}{cc} D^2 \phi(\xe)
& 0 \\ 0 & 0
\end{array}\right) ,$$
and
$$ F(\e \Xt,\pt,\frac{\xe}{\e}, \xe)+\ue(\xe)\leq 0 \; ,$$
and
$$ F_k (\Yt, \qt + D\phi(\xb) , \frac{\ze}{\e}, \xb) \geq \Fb_k (D \phi(\xb),\xb)\; .$$
Moreover $\pt = D\phi(\xe) + \qt$. From the matrix inequality, we deduce that $\displaystyle \Xt \leq \frac 1 \e \Yt + D^2 \phi(\xe)$ and by standard arguments, 
$\displaystyle \frac{\vert \xe -\ze \vert^2}{\beta^2}$ is uniformly bounded. Using these properties and the fact that $\beta \ll \e$, we have
\begin{eqnarray*}
 F_k (\Yt, \qt + D\phi(\xb) , \frac{\ze}{\e}, \xb) & \leq & F_k (\e \Xt-\e D^2 \phi(\xe) , \pt - D\phi(\xe) + D\phi(\xb) , \frac{\ze}{\e}, \xb)\\
 & \leq & F_k (\e \Xt +o(1) , \pt + o(1) , \frac{\xe}{\e}+o(1), \xe + o(1))\\
 & \leq & F (\e \Xt , \pt  , \frac{\xe}{\e}, \xe )
\end{eqnarray*}
where the last inequality comes from the definition of $F_k$ as soon as the $o(1)$ are less than $k^{-1}$.

Therefore
\begin{eqnarray*}
 \Fb_k (D \phi(\xb),\xb) + \ue(\xe) & \leq & F_k (\Yt, \qt + D\phi(\xb) , \frac{\ze}{\e}, \xb) + \ue(\xe) \\
& \leq & F (\e \Xt , \pt  , \frac{\xe}{\e}, \xe ) +\ue(\xe) \\
& \leq & 0  \; ,
\end{eqnarray*}
and by sending $\e$ to $0$, we obtain $ \Fb_k (D \phi(\xb),\xb) + \ou(\xb) \leq 0$. And we conclude by letting $k$ tends to infinity, since by Theorem~\ref{erg-pb}, we know that $\Fb_k (D \phi(\xb),\xb) \to \Fb (D \phi(\xb),\xb)$.

\section{Appendix}

The aim of this appendix is to study the existence of the functions $\chi_1, \chi_2$ used in the proof of Theorem~\ref{mainregu}.

Before doing that, we remark that, if $h$ is in the class $\mathcal{P}$ and if $t_2 \geq t_1 >0$, then the property (ii) of the class $\mathcal{P}$ implies
$$ \frac{h(t_2)}{t_2^2} \geq \frac{h(t_1)}{t_1^2}\; ,$$
and by taking (for example) $t_1 = 1$ and $t_2$ large, we have
$$ h(t_2) \geq t_2^2 h(1)\; ,$$
which means that $h$ grows at least quadratically at infinity. On the other hand, choosing $t_2 = Ct$, $t_1=t$ for some $C \geq 1$, $t>0$, we get
$$\frac{h(Ct)}{C^2 t^2} \geq \frac{h(t)}{t^2}\quad\hbox{i.e.}\quad h(Ct) \geq C^2 h(t)\; .$$

Then we first examine the ode
$$\chi''_2 (t) = - h(\chi '_2 (t))\; .$$
Using the fact that $ \chi'_2(0)=+\infty$, we deduce that, necessarely, $\chi '_2 (t)$ solves
$$ \int_{\chi '_2 (t)}^{+\infty}\,\frac{ds}{h(s)} = t\; .$$
This equation has the form $F(\chi '_2 (t)) = t$ where $F$ is the strictly decreasing function
$$ F(\tau) = \int_{\tau}^{+\infty}\,\frac{ds}{h(s)}\; .$$
Solving the ode just consists in inverting this function (changing a little bit $h$ if necessary in order to be able to do it for $t$ in $[0,1]$ or $[0,2]$).

To check that we can choose $\chi_2(0)=0$, we come back to the original ode which can be written as
$$ \frac{\chi '_2 (t)}{ h(\chi '_2 (t))}\chi''_2 (t) = -\chi '_2 (t)\; ,$$
and integrating from $\e \ll 1$ to $1$, yields
$$ \chi _2 (\e)-\chi_2 (1) = \int_{ \chi _2 (\e)}^{ \chi _2 (1)}\frac{s}{h(s)}ds\; .$$
The integral in the right-hand side being convergent, $\chi _2 (\e)$ has a limit as $\e\to 0$ and adjusting the constant we can assume that  $\chi_2 (0)=0$.  

For $\chi_1$, we are in the opposite situation since we start from what is above the equality $F(\chi '_2 (t)) = t$ for $\chi_2$. We set $\displaystyle g(\tau) = \frac{\tau}{h(\tau)}$. The equation 
$$ g(\chi'_1 (t)) = t\; ,$$
can be solved since by the assumption (ii) of the class $\mathcal{P}$,  $\displaystyle t \mapsto \frac{h(t)}{t^2}$ is non-decreasing for $t\geq 1$ and this is a fortiori the case for  $\displaystyle t \mapsto \frac{h(t)}{t}$. It is clear that $\chi'_1 (0)=+\infty$ and for $\chi_1 (0)$, we repeat the above arguments : by differentiating the $\chi'_1$-equation and multiplying by $\chi'_1 (t)$, we have
$$ \chi''_1 (t)g'(\chi'_1 (t))\chi'_1 (t) = \chi'_1 (t)\; ,$$
and we follow the same step as for $\chi_2$. We point out that
$$ g'(\tau)\tau = \frac{\tau}{h(\tau)} - \frac{\tau^2 h'(\tau)}{[h(\tau)]}\; ,$$
which has the right integrability property at infinity.

\end{document}